\documentclass[12pt]{article}   

\usepackage{scicite}
\usepackage{times}
\usepackage{subcaption}

\usepackage{graphicx}				% Use pdf, png, jpg, or eps§ with pdflatex; use eps in DVI mode
\usepackage{amssymb}
\usepackage{multirow}
\usepackage{comment}
\usepackage{amsmath}
\newcommand{\pvalue}{{\em P}-value}
\newcommand{\pvalues}{{\em P}-values}
\newcommand{\bt}{{\boldsymbol \theta}}
\newcommand{\bfbeta}{{\boldsymbol \beta}}
\newcommand{\bfomega}{{\boldsymbol \omega}}
\newcommand{\bfo}{{\boldsymbol \omega}}
\newcommand{\bP}{{\bold P}}

\newcommand{\con}{{\,|\,}}

\newcommand{\bfx}{{\bold x}}

\newcommand{\bp}{{\bf p}}

\newcommand{\kt}{{\frac{k}{2}}}

\newcommand{\bpi}{{\boldsymbol \pi}}

\newenvironment{sciabstract}{%
\begin{quote} \bf}
{\end{quote}}

\title{Bayes factor functions for reporting outcomes of hypothesis tests}

\author
{Valen E. Johnson,$^{1\ast}$ Sandipan Pramanik,$^{1}$ Rachael Shudde$^{1}$\\
\\
\normalsize{$^{1}$Department of Statistics, Texas A\&M University,}\\
\normalsize{3143 TAMU, College Station, TX 77843-3143, USA}\\
\\
\normalsize{$^\ast$To whom correspondence should be addressed; E-mail:  vejohnson@tamu.edu}
}

\date{}

\begin{document}

\maketitle

\begin{sciabstract}
Bayes factors represent the ratio of probabilities assigned to data by competing scientific hypotheses.  Drawbacks of Bayes factors are their dependence on prior specifications that define null and alternative hypotheses and difficulties encountered in their computation.  To address these problems, we define Bayes factor functions (BFF) directly from common test statistics.  BFFs depend on a single non-centrality parameter that can be expressed as a function of standardized effect sizes, and plots of BFFs versus effect size provide informative summaries of hypothesis tests that can be easily aggregated across studies. Such summaries eliminate the need for arbitrary \pvalue\ thresholds to define ``statistical significance.''  BFFs are available in closed form and can be computed easily from $z$, $t$, $\chi^2$, and $F$ statistics.  
\end{sciabstract}

Two approaches are commonly used to summarize evidence from statistical hypothesis tests: \pvalues\  and Bayes factors.  $P$-values are more frequently reported.  As noted in the American Statistical Association Statement on Statistical Significance and {\em P}-values, the significance of many published scientific findings are based on \pvalues, even though this index ``is commonly misused and misinterpreted.  This has led to some scientific journals discouraging the use of $P$-values, and some scientists and statisticians recommending their abandonment ... Informally, a $P$-value is the probability under a specified statistical model that a statistical summary of the data would be equal to or more extreme than its observed value." \cite{ASA}.   \pvalues\  do not provide a direct measure of support for either the null or alternative hypotheses, and their use in defining arbitrary thresholds for defining statistical significance has long been a subject of intense debate; see, for example, \cite{Edwards63, Berger87,Johnson2013pnas, Nuzzo2014, Greenland2016, Benjamin2017,Lakens2018}.  Interpreting evidence provided by \pvalues\ across replicated studies can also be challenging.

Bayes factors represent the ratio of the marginal probability assigned to data by competing hypotheses and, when combined with prior odds assigned between hypotheses, yield an estimate of the posterior odds that each hypothesis is true. That is, 
\begin{equation}\label{bayes1}
\mbox{posterior odds}  = \mbox{Bayes factor}  \times   \mbox{prior odds},    
\end{equation}
or, more precisely,
\begin{equation}\label{bayes2}
\frac{\bP(H_1\con \bfx)}{\bP(H_0 \con \bfx)}  =   \frac{m_1(\bfx)}{m_0(\bfx)}  \times  \frac{\bP(H_1)}{\bP(H_0)}.
\end{equation}
Here, $\bP(H_i \con \bfx)$ denotes the posterior probability of hypothesis $H_i$ given data $\bfx$; $\bP(H_i)$ denotes the prior probability assigned to $H_i$; and $m_i(\bfx)$ denotes the marginal probability (or probability density function) assigned to the data under hypothesis $H_i$, for $i=0$ (null) or $i=1$ (alternative).

The marginal density of the data under the alternative hypothesis is given by
\begin{equation}
m_1(\bfx) = \int_\Theta f(\bfx \con \theta) \pi_1(\theta) d\,\theta.
\end{equation}
In null hypothesis significance tests (NHSTs), the marginal density of the data under the null hypothesis, $m_0(\bfx)$, is simply the sampling density of the data assumed under the null hypothesis.
The function $\pi_1(\theta)$ represents the prior density for the parameter of interest $\theta$ under the alternative hypothesis, i.e., the alternative prior density.  

The specification of $\pi_1(\theta)$ is problematic and, as a consequence, numerous Bayes factors based on ``default'' alternative prior densities have been proposed.  These include \cite{Ohagan95, Berger1996, Berger96,Liang2008,Rouder2009,Wagenmakers2010,Consonni2018,Etz2018}.  Nonetheless, the value of a Bayes factor depends on the alternative prior density used in its definition, and it is generally difficult to justify or interpret any single default choice. In addition, the numerical calculation of many Bayes factors is difficult, often requiring specialized software, and each of these problems is exacerbated in high-dimensional settings \cite{Morey2011}.  
A more detailed description of the Bayesian hypothesis testing framework and controversies surrounding the definition of default alternative prior densities may be found in, for example, \cite{Jeffreys1961} or \cite{Kass1995}.  

We propose several modifications to existing Bayes factor methodology designed to enhance the report of scientific findings.  First, we define Bayes factors directly from standard $z$, $t$, $\chi^2$, and $F$ test statistics \cite{Johnson2005}.  Under the null hypotheses, the distribution of these test statistics is known.  Under alternative hypotheses, the asymptotic distributions of these test statistics often depend only on scalar-valued non-centrality parameters.  Thus, the specification of the prior density that defines the alternative hypothesis is simplified, and no prior densities need to be specified under the null hypothesis.  

Second, we express prior densities on non-centrality parameters as functions of standardized effect sizes and construct  Bayes Factor Functions (BFFs) that vary with these standardized effects.  BFFs thus make the connection between Bayes factors and prior assumptions more apparent.   BFFs also facilitate the integration of evidence across multiple studies of the same phenomenon.  

Third, the prior densities we propose for non-centrality parameters are special cases of non-local alternative prior (NAP) densities.  These densities are identically zero when the non-centrality parameter is zero.  This property makes it possible to more quickly accumulate evidence in favor of both true null and true alternative hypotheses \cite{Johnson2010,Rossell2017,Cao21}; this particular feature of NAP densities is discussed in detail in \cite{Pramanik2023}.

Finally, we provide closed-form expressions for Bayes factors and BFFs.  These expressions eliminate computational difficulties sometimes encountered when calculating other Bayes factors.  

 \section{Mathematical framework}
Theorems describing default Bayes factors based on $z$, $t$, $F$, and $\chi^2$ statistics are provided below.  In each case, the Bayes factors depend on a hyperparameter $\tau^2$.   Procedures for setting $\tau^2$ as a function of standardized effect size are described in Section~\ref{criteria}.  Proofs of the theorems appear in the Supplemental Material.

Notationally, we write $a \con {\bf b} \sim D({\bf b})$ to indicate that a random variable $a$ has distribution $D$ that depends on a parameter vector ${\bf b}$.  $N(a,b)$ denotes the normal distribution with mean $a$ and variance $b$; $T(\nu,\lambda)$ denotes a $t$ distribution with $\nu$ degrees-of-freedom and non-centrality parameter $\lambda$; $F(k,m,\lambda)$ denotes an $F$ distribution on $k,m$ degrees-of-freedom and non-centrality parameter $\lambda$; $G(\alpha,\lambda)$ denotes a gamma random variable with shape parameter $\alpha$ and rate parameter $\lambda$; $H(\nu,\lambda)$ denotes a $\chi^2$ distribution on $\nu$ degrees of freedom and non-centrality parameter $\lambda$; and $J(\mu_0,\tau^2)$ denotes a normal moment distribution with mean $\mu_0$ and rate parameter $\tau^2$ \cite{Johnson2010}.
We use lower case letters to denote corresponding densities; e.g., a gamma density evaluated at $x$ is written $g(x\con \alpha, \lambda)$.  

\begin{figure}[ht]
\begin{subfigure}{.5\textwidth}
  \centering
  % include first image
  \includegraphics[width=.8\linewidth]{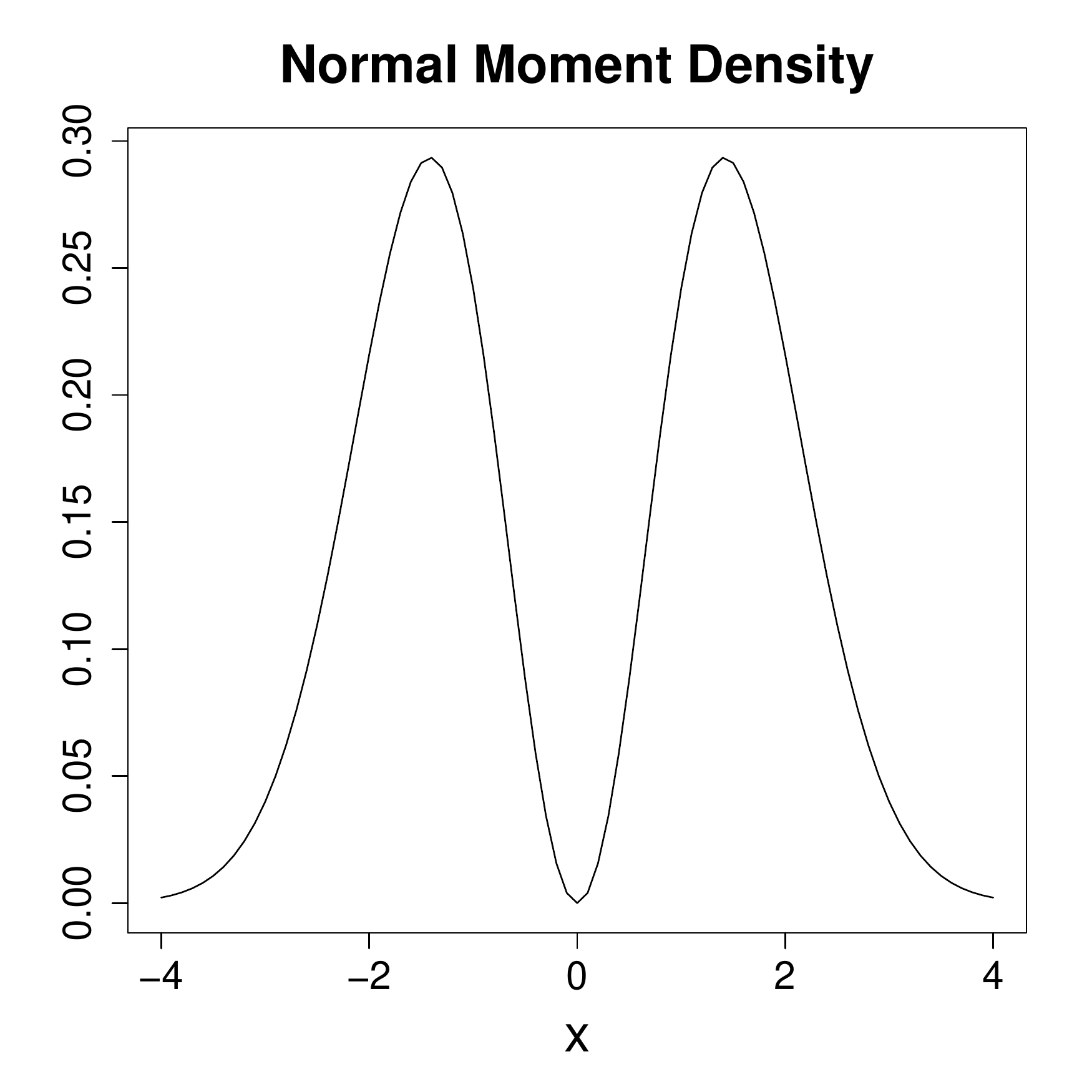} 
\end{subfigure}
\begin{subfigure}{.5\textwidth}
  \centering
  % include second image
  \includegraphics[width=.8\linewidth]{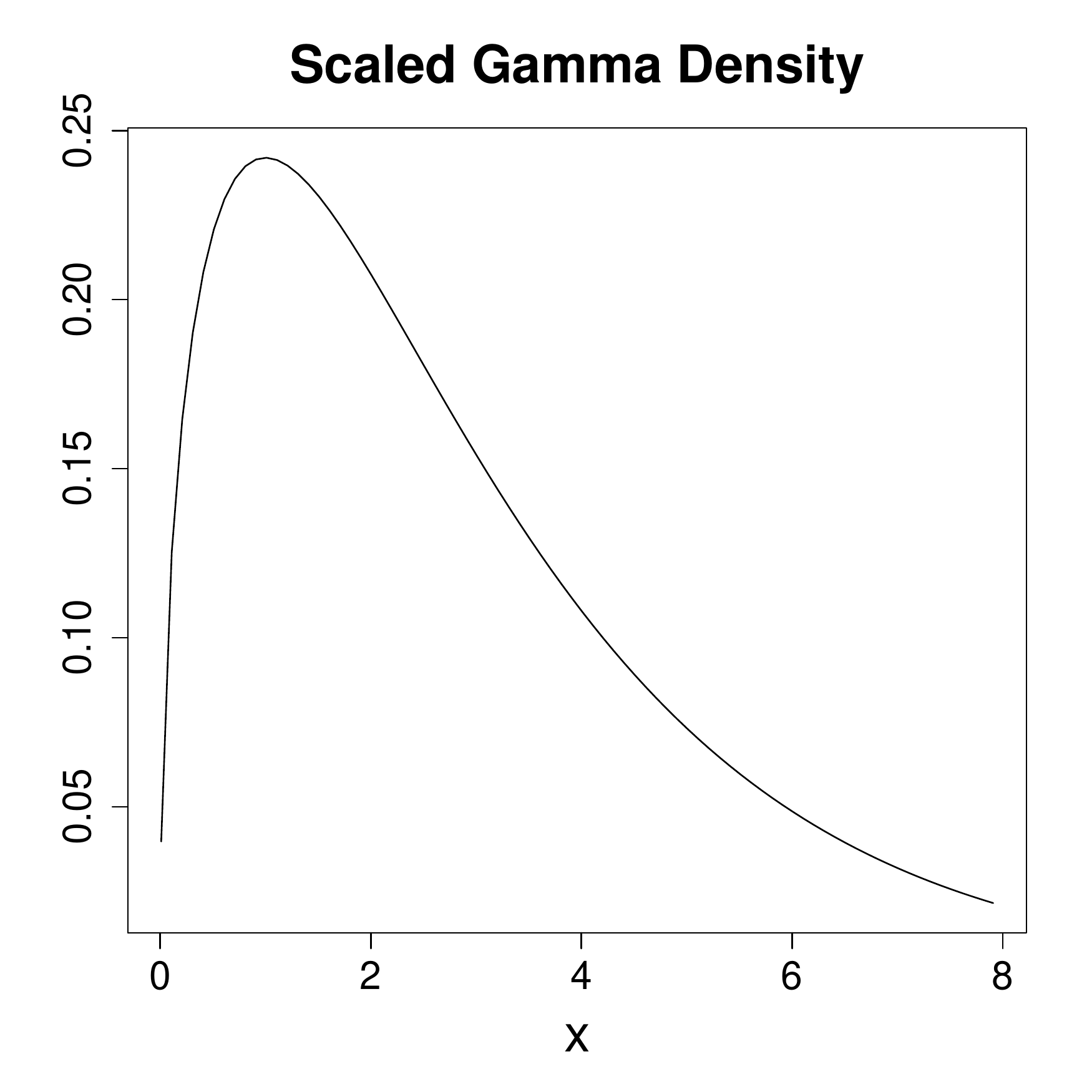}  
\end{subfigure}
\caption{Normal moment prior density, $J(0,1)$ (left panel), and scaled Gamma density, $G\left(\frac{3}{2}, \frac{1}{2}\right)$ (right panel).}
\label{napplot}
\end{figure}

The density of a $J(\mu_0,\tau^2)$ random variable can be written as 
\begin{equation}
j( x \con \mu_0, \tau^2 ) = \frac{(x-\mu_0)^2}{\sqrt{2 \pi} \tau^3  } \exp\left( -\frac{(x-\mu_0)^2}{2\tau^2} \right).
\end{equation}
The modes of this distribution occur at $x=\mu_0 \pm \sqrt{2} \tau$.   A plot of a normal moment density appears in the left panel of Figure~\ref{napplot}.  If $a\sim J(0,\tau^2)$, then $a^2$ has a $G[3/2,1/(2\tau^2)]$ distribution, as displayed in the right panel of Figure~\ref{napplot}.

The first two theorems extend results in \cite{Pramanik2023} for parametric hypothesis tests for normally distributed data to more general cases of $z$ and $t$ tests.  Theorems~3 and~4 provide new results for $\chi^2$ and $F$ tests. Throughout, $BF_{10}(x\con \psi)$ denotes the Bayes factor in favor of the alternative hypothesis based on a test statistic $x$ for a hyperparameter value $\psi$.

{\noindent \bf Theorem 1: z test}. {\em Assume the distributions of a random variable $z$ under the null and alternative hypotheses are described by
\begin{eqnarray}
H_0: z &\sim& N(0,1), \\
H_1: z\con \lambda &\sim& N(\lambda,1), \qquad \lambda\con \tau^2 \sim J(0,\tau^2).
 \end{eqnarray}
 Then the Bayes factor in favor of the alternative hypothesis is
 \begin{equation}\label{zBF}
 BF_{10}(z\con \tau^2) = (\tau^2+1)^{-3/2} \left(1+\frac{\tau^2 z^2}{\tau^2+1}  \right) \exp\left[\frac{\tau^2 z^2}{2(\tau^2+1)}\right].
 \end{equation}
 }

{\noindent \bf Theorem 2: t test}. {\em  Assume the distributions of a random variable $t$ under the null and alternative hypotheses are described by
\begin{eqnarray}
H_0: t &\sim& T(\nu,0), \\
H_1: t \con \lambda &\sim&
 T(\nu,\lambda), \qquad \lambda \con \tau^2 \sim J(0,\tau^2).
 \end{eqnarray}
 Then the Bayes factor in favor of the alternative hypothesis is
 \begin{equation}\label{tBF}
 BF_{10}(t\con \tau^2) = (\tau^2+1)^{-3/2} \left( \frac{r}{s} \right)^{\frac{\nu+1}{2}}  \left( 1 +\frac{qt^2}{s} \right),
 \end{equation}
 where
 \[ r = 1 + \frac{t^2}{\nu}, \qquad s = 1+ \frac{t^2}{\nu(1+\tau^2)}, \quad \mbox{and} \quad q= \frac{\tau^2(\nu+1)}{\nu(1+\tau^2)}. \]}

{\noindent \bf Theorem 3: $\chi^2$ test}.  {\em Assume the distributions of a random variable $h$ under the null and alternative hypotheses are described by 
\begin{eqnarray}
H_0:  h &\sim& H(k,0), \\
H_1: h \con \lambda &\sim& H(k,\lambda), \qquad \lambda \con \tau^2 \sim G\left( \frac{k}{2} + 1, \frac{1}{2\tau^2} \right). \label{chisqprior}
 \end{eqnarray}
Then the Bayes factor in favor of the alternative hypothesis is
\begin{equation}\label{x2BF}
BF_{10}(h \con \tau^2)  = (\tau^2+1)^{-k/2-1}  \left( 1 + \frac{\tau^2}{k(\tau^2+1)} h \right) \exp\left( \frac{\tau^2 h}{2(\tau^2+1)} \right).
\end{equation}
}

For $k=1$ and $z^2=h$, the Bayes factors in equation (\ref{x2BF}) has the same value as the Bayes factor specified in equation (\ref{zBF}).  The choice of the shape parameter as $k/2+1$ for the gamma density (a scaled $\chi^2_{k+2}$ random variable) in equation (\ref{chisqprior}) was based on the fact that $\chi^2_{\nu}$ distributions are not 0 at the origin for integer degrees of freedom $\nu<3$.  Thus, they are not NAP densities for $\nu<3$ and typically are not able to provide strong evidence in favor of true null hypotheses without large sample sizes \cite{Johnson2010}.

{\noindent \bf Theorem 4: F test}.  {\em Assume the distributions of a random variable $f$ under the null and alternative hypotheses are described by
\begin{eqnarray}
H_0:  f &\sim& F(k,m,0), \\
H_1: f \con \lambda &\sim& F(k,m,\lambda), \qquad \lambda \con \tau^2 \sim G\left( \frac{k}{2} + 1, \frac{1}{2\tau^2} \right).\label{falt}
 \end{eqnarray}
Then the Bayes factor in favor of the alternative hypothesis is
\begin{equation}\label{fBF}
BF_{10}(f \con \tau^2)  = (\tau^2+1)^{-\frac{k}{2}-1}  
\left[ \frac{\left(1+\frac{kf}{m} \right) }{ \left(1+\frac{kf}{v} \right) } \right]^{\frac{k+m}{2}} 
\left[ 1 + \frac{\, (k+m)\tau^2 \,f}{v\, \,\left(1 + \frac{kf}{v} \right)} \right],
\end{equation}
where $v=m(\tau^2+1)$.}

For $k=1$ and $t^2=f$, the Bayes factors in equation (\ref{fBF}) has the same value as the Bayes factor specified in equation (\ref{tBF}).

\subsection{Bayes factors as functions of standardized effect size}\label{criteria}
Theorems 1-4 describe Bayes factors based on classical test statistics.  
In each case, the distribution of the test statistic under the alternative hypothesis depends on a prior distribution on a non-centrality parameter, which in turn depends on a hyperparameter $\tau^2$.  The non-centrality parameters determine the deviation of the alternative hypothesis from the null hypothesis. Like the hyperparameters that define priors in other Bayesian tests, the value of $\tau^2$ plays an essential role in determining the Bayes factor. Rather than ignoring this dependence or simply reporting a single Bayes factor that depends on a specific hyperparameter value, we construct BFFs that vary with $\tau^2$.  Unfortunately, $\tau^2$ is difficult to interpret scientifically and, as we show below, its interpretation changes from one test to another.  For this reason, we report BFFs as functions of standardized effect sizes. 

In typical applications of $z$, $t$, $\chi^2$, and $F$ tests, standardized effect sizes determine the non-centrality parameters of the test.  For example, consider 
a $z$ test of a null hypothesis $H_0: \mu=0$ based on a random sample $x_1,\dots,x_n$, where $x_i \sim N(\mu,\sigma)$  and $\sigma^2$ is known. The $z$ statistic for this test is $z=\sqrt{n} \bar{x}/\sigma$.  The distribution of the $z$ statistic is
\begin{equation}\label{zex}
z\con \mu, \sigma \sim N\left(\frac{\sqrt{n}\mu}{\sigma},1\right).
\end{equation}
Under the null hypothesis, $z \sim N(0,1)$.  Under the alternative hypothesis, $\mu$ defines the deviation from the null value 0 and is called the effect size. The non-centrality parameter  is $\lambda= \sqrt{n}\mu/\sigma$.  Thus, the non-centrality parameter varies with both $n$ and $\mu/\sigma$. 

In many studies, effect sizes are standardized.  Standardization serves two purposes.  First, it makes the effect size invariant to units of measurement--for example, whether weights are measured in ounces or grams.  Second, standardization scales the effect size according to the random variation between observational units.  For the $z$ test, we can standardize the effect size $\mu$ by dividing it by the standard deviation of the observations, leading to a standardized effect $\omega = \mu/\sigma$.  The non-centrality parameter for the test can then be expressed as $\lambda = \sqrt{n} \omega$ \cite{Cohen88,meta2021}.

Given a relationship between a standardized effect size $\omega$ and a non-centrality parameter $\lambda$, we compute the BFF by setting $\tau^2$ so that the modes of the prior density on $\lambda$ occur at values defined by standardized effect sizes $\omega$.

More explicitly, suppose that the non-centrality parameter $\lambda$ can be written as a function of the standardized effect size $\omega$ as $\lambda = r(\omega)$, and let $\pi(\lambda\con \tau^2)$ denote a generic prior density on $\lambda$ given $\tau^2$.  Define $\tau^2_\omega$ implicitly as 
\begin{equation}
r(\omega) = \underset{\lambda}{\arg\max} \ \pi(\lambda \con \tau^2_\omega), 
\end{equation}
the value of $\tau^2$ that makes the prior modes equal to $r(\omega)$. 
Given $\tau^2_\omega$ for a range of $\omega$ values, the BFF based on $x$ consists of ordered pairs $(BF(x \con \tau^2_\omega),\omega)$.

To illustrate this procedure, consider again the test of a normal mean. 
The non-centrality parameter for this test is $\lambda = \sqrt{n}\omega$.  The default prior on the non-centrality parameter $\lambda$ is a $J(0,\tau^2)$ distribution, which has maxima at $\lambda = \pm \sqrt{2}\tau$ (that is, $\pm \sqrt{2}\tau = \underset{\lambda}{\arg\max} \ j(\lambda\con 0,\tau^2)$).  
Equating the non-centrality parameter to the modes of the prior density (i.e., $\sqrt{n} \omega = \lambda = \pm \sqrt{2}\tau$) implies that $\tau_\omega^2 = n\omega^2/2$.

\begin{figure}[htbp] 
   \centering
   \includegraphics[width=3.5in,height=3in]{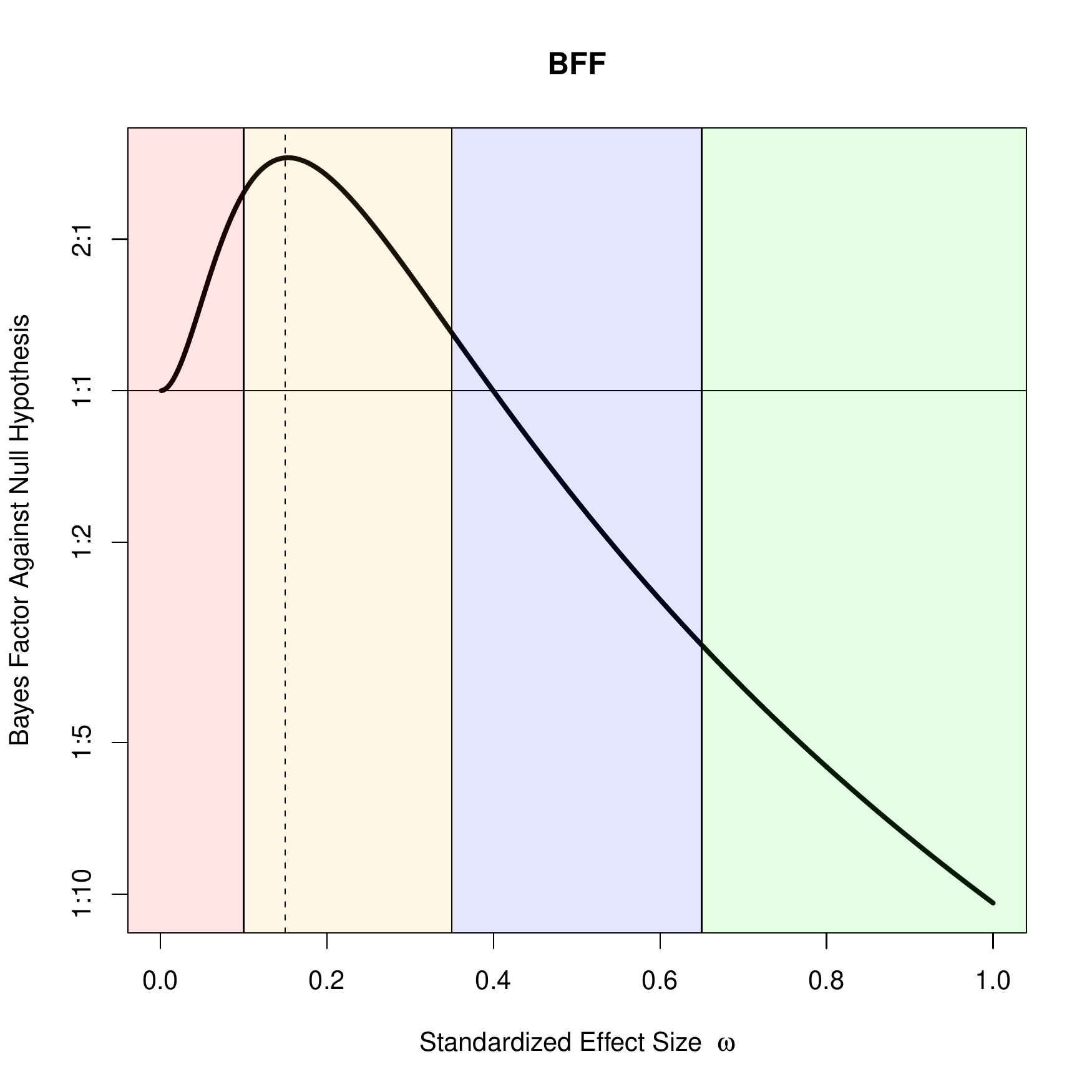} 
   \caption{Plot of the BFF, $BF_{10}(2 \con 100\omega^2/2)$, against $\omega$ for a $z$ test with $z=2$ and $n=100$. Bayes factors are displayed as odds in favor of the alternative hypothesis.  The vertical axis is displayed on the logarithmic scale.  The vertical dotted line indicates that the maximum Bayes factor of 2.90 corresponds to a standardized effect size of 0.15.  }
   \label{fig:one}
\end{figure}

Figure~\ref{fig:one} displays the BFF in favor of the alternative hypothesis using the mapping $\tau_\omega^2=n\omega^2/2$ when $z=2.0$ and $n=100$.  From the BFF, we can conclude that the maximum Bayes factor in favor of the alternative hypothesis is 2.90 and that this Bayes factor is achieved when the prior modes on the standardized effect size are $\pm 0.15$.  Furthermore, the odds in favor of the alternative hypothesis fall below 1:1 for standardized effect sizes greater than $0.4$ and are below 1:5 when the prior modes of the distribution on effect size are greater than 0.8.  As the standardized effect size approaches 0, the alternative hypothesis becomes indistinguishable from the null hypothesis, driving the Bayes factor to 1.  The red, orange, blue, and green zones in this figure are arbitrarily colored and correspond to very small ($(0,0.1)$), small ($(0.1,0.35)$), medium ($(0.35,0.65)$) and large ($>0.65$) standardized effect sizes, respectively.

Table~\ref{tauTable} provides a mapping between standardized effect sizes $\omega$ and $\tau^2_\omega$ for several common statistical tests.   Special cases of the tests in the ``Multinomial/Counts" row include Pearson's $\chi^2$ test for goodness-of-fit ($s=0$ and $f(\bt)$ known) and tests for independence in contingency tables ($K-s-1 = (\mbox{\# rows}-1)(\mbox{\# columns}-1)$).  Recall that a test for the value of a binomial proportion can be framed as a Pearson's $\chi^2$ test.  In contrast, a test for a difference in proportions can be specified as a test for independence in contingency tables.  

The last three rows in Table 1 contain vectors of standardized effects $\bfo$.  Because the recommended value of $\tau^2_\omega$ depends on $\bfo$ only through the inner product $\bfo'\bfo$, in many applications it is easier to study the BFF as a function of the root mean square effect size (RMSES), $\tilde{\omega}$, defined as

\begin{equation}
\tilde{\omega} = \sqrt{\frac{1}{k} \sum_{i=1}^k \omega_i^2}.
\end{equation}

The last entry in Table~1 provides default choices for $\tau_\omega^2$ based on the asymptotic distribution of the likelihood ratio statistic.  This choice is based on classical results summarized in, for example, \cite{Kendall1991}.  This entry is of particular interest due to the widespread application of the likelihood ratio test statistic in non-linear models.  

Justification for the values of $\tau_\omega^2$ in Table~1\footnote{Notation in Table~1.  For one-sample tests, $x_1,\dots,x_n$ are assumed to be iid $N(\mu,\sigma^2)$ where $n$ refers to sample size.  In two-sample tests, $x_{j,1},\dots,x_{j,n_j}$, $j=1,2$, are assumed to be iid $N(\mu_j,\sigma^2)$. Integers $n_1$ and $n_2$ refer to sample sizes in each group. A bar over a variable denotes the sample mean. The variance of normal observations is denoted by $\sigma^2$ and is assumed to be equal in both groups in two-sample tests. Standard deviations are denoted by $s$ and are the pooled estimate in the two-sample $t$ test. In Multinomial/Poisson tests, $f(\bt)$ maps an $s \times 1$ vector $\bt$ into a $k \times 1$ probability vector, where $k$ denotes the number of cells. The degrees of freedom $\nu$ equals $k-s-1$.  The quantities $p_i$ and $n_i$ represent cell probabilities and counts, respectively, and $n$ is the sum of all cell counts.  In the Linear Model, the alternative hypothesis is ${\bf A}\bfbeta = {\bf a}$, where {\bf A} is a $k\times p$ matrix of rank $k$, $\bfbeta$ is a $p\times 1$ vector of regression coefficients, and {\bf a} is a $k\times 1$ vector. The quantities $RSS_0$ and $RSS_1$ denote the residual sum-of-squares under the null and alternative hypotheses, respectively.  The quantity $n$ is the number of observations, and $\sigma^2$ is the observational variance. In the Likelihood Ratio test, $l(\cdot)$ denotes the likelihood function for a parameter vector $\bt = (\bt_r,\bt_s)$. The $k\times 1$ subvector $\bt_r$ equals $\bt_{r0}$ under the null hypothesis. The maximum likelihood estimate of $\bt$ under the alternative hypothesis is $\hat{\bt}$, and the maximum likelihood of $\bt_s$ under the null hypothesis is $\hat{\bt}_s$.  In the Linear Model and Likelihood Ratio tests, the matrix ${\bf L}^{-1}$ represents the Cholesky decomposition of the covariance matrix for the tested parameters, scaled to a single observation.  Further explanation of $\tau^2_{\omega}$ values appear in the Supplemental Material.} appears in the Supplemental Material.

\begin{table}[h]
\def\arraystretch{2.0}
{
\begin{tabular}{lccc} \hline \hline
{\bf Test }& Statistic  & Standardized Effect ($\omega$) & $\tau^2_\omega$     \\ \hline \hline
{1-sample z} & {$\frac{\sqrt{n}\bar{x}}{\sigma}$} &  $\frac{\mu}{\sigma}$ & $ \frac{n\omega^2}{2}$
           \\  \hline
{1-sample t} & {$\frac{\sqrt{n}\bar{x}}{s}$} &  $\frac{\mu}{\sigma}$ & $ \frac{n\omega^2}{2}$
           \\  \hline
{2-sample z} 
	& $\frac{\sqrt{n_1 n_2}(\bar{x}_1-\bar{x}_2)}{\sigma\sqrt{n_1+n_2}} $ & $\frac{\mu_1-\mu_2}{\sigma}$ 
         & $\frac{n_1 n_2\omega^2}{2(n_1+n_2)}$ \\ 
\hline
{2-sample t} 
	& $\frac{\sqrt{n_1 n_2}(\bar{x}_1-\bar{x}_2)}{s\sqrt{n_1+n_2}} $ & $\frac{\mu_1-\mu_2}{\sigma}$  
         & $\frac{n_1 n_2\omega^2}{2(n_1+n_2)}$ \\ 
\hline

Multinomial/Poisson & 
          \multirow{2}{*}{$\chi^2_{\nu}= \sum\limits_{i=1}^k \frac{(n_i-nf_i(\hat{\bt}))^2}{nf_i(\hat{\bt})}$} &
          \multirow{2}{*}{{ $ \left( \frac{p_{i}-f_i(\bt)}{\sqrt{f_i(\bt)}} \right)_{k\times 1} $} }& 
          \multirow{2}{*}{$ \frac{n \bfo'\bfo}{k} = n\tilde{\omega}^2$} \\
 & & &   \\ \hline
Linear model& \multirow{2}{*}{$F_{k,n-p} = \frac{(RSS_0-RSS_1)/k}{[(RSS_1)/(n-p)]}$} &
           \multirow{2}{*}{ $ \frac{{\bf L}^{-1}({\bf A}\bfbeta-{\bf a})}{\sigma} $} & 
           \multirow{2}{*}{ $\frac{n\bfo'\bfo}{2k} = \frac{n\tilde{\omega}^2}{2} $}
                \\
                \\ \hline
Likelihood Ratio & 
       \multirow{2}{*}{$ \chi^2_{k} = -2\log\left[\frac{l(\bt_{r0},\hat{\bt_{s}})}{l(\hat{\bt})} \right]$} &
       \multirow{2}{*}{${\bf L}^{-1}(\bt_{r}-\bt_{r0}) $ } & 
       $\frac{n{\bfomega'}{\bfomega}}{k}=n\tilde{\omega}^2$ \\
 &  \\  
     \hline
\end{tabular}
}
\caption{Default choices for $\tau^2_\omega$. 
}\label{tauTable}
\end{table}

\section{Applications}
The following examples show how BFF can be used to summarize outcomes of hypothesis tests based on $\chi^2$ and $F$ test statistics.
\subsection{Cancer sites and blood type association}
White and Eisenberg \cite{White1959} collected data from 707 patients with stomach cancer and investigated the association between cancer site and blood type. Data from their study are summarized in Table 2.  The $\chi^2$ test for independence for these data is 12.65 on 6 degrees of freedom.
\begin{table}
\begin{center}
\begin{tabular}{lccc}\hline
Site & \multicolumn{3}{c}{Results for the} \\
 & \multicolumn{3}{c}{following blood groups}  \\ \hline
 & O & A & B or AB \\ \hline
 Pylorus and antrum & 104 & 140 & 52 \\
 Body and fundus & 116 & 117 & 52  \\
 Cardia & 28 & 39 & 11 \\
 Extensive & 28 & 12 & 8 \\ \hline
 \end{tabular}
 \caption{White and Eisenberg's classification of cancer patients}
 \end{center}
 \end{table}
Bayes factors to test the independence of cancer site and blood type were previously calculated by Albert \cite{Albert1990}, Good and Crook \cite{Good1987}, and Johnson \cite{Johnson2005}.  The Bayes factors reported in \cite{Albert1990,Good1987} require the specification of prior distributions on the marginal probabilities of blood type and cancer site under the null hypothesis and the specification of a Dirichlet distribution on all combinations of $\mbox{blood type}\times\mbox{cancer site}$ probabilities under an alternative model. As noted in \cite{Johnson2005}, the specification of these prior probability models is ``rather intricate;'' interested readers should consult the original articles for details.  Johnson maximized a Bayes factor based on the chi-squared statistic similar to that proposed in Theorem~3, except that the prior on the non-centrality parameter was a scaled chi-squared distribution on $6$ degrees of freedom (rather than 8).  The scale of the chi-squared prior was chosen to maximize the Bayes factor against the null hypothesis of independence.  The Bayes factors reported in \cite{Johnson2005,Albert1990,Good1987} were 2.97, 3.02, and 3.06, respectively.

\begin{figure}[ht] 
   \centering
   \includegraphics[width=3.5in]{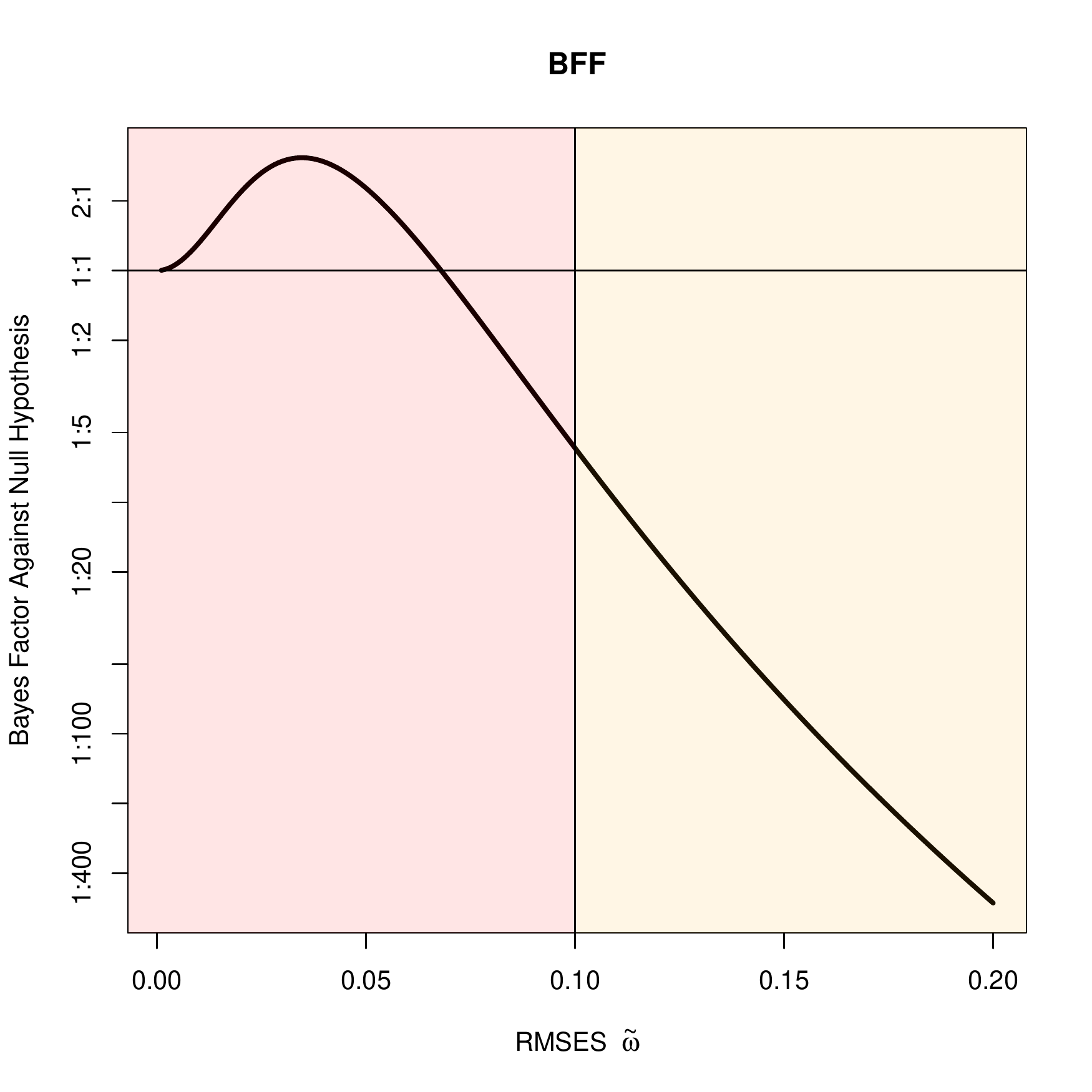} 
   \caption{Plot of the BFF, $BF_{10}(12.65 \con 707\tilde{\omega}^2)$, against $\tilde{\omega}$ for a $\chi^2_6$ test with $h=12.65$.  Bayes factors are displayed as odds in favor of dependence between patient blood type and cancer site.  The vertical axis is displayed on the logarithmic scale.  The color coding is consistent with Fig.~2.  The horizontal line in the plot corresponds to a Bayes factor of 1.0 (odds of 1:1).  }
   \label{fig:chi2}
\end{figure}

Fig.~\ref{fig:chi2} displays the BFF as a function of the RMSES using results from Theorem~3 and the $\tau^2_\omega$ values provided in Table~1.  The maximum Bayes factor in favor of dependence occurs for $\tilde{\omega}=0.035$, where it equals 3.07.   The Bayes factor favors the independence model against alternatives with $\tilde{\omega}>0.068$.  A standardized effect size of 0.2 represents what is often considered a small standardized effect in the social science and medical literature \cite{Cohen88}, and the Bayes factor against such an effect size for these data is greater than 400:1.  
The Bayes factors reported in \cite{Johnson2005,Albert1990,Good1987} and the maximum value of 3.07 reported here are, for many practical purposes, similar in their scientific interpretation.  They all suggest that the data support an alternative hypothesis of non-independence three times more than the null hypothesis of independence.   However, the previous methods do not provide a straightforward interpretation of the alternative hypothesis they support in terms of standardized effect sizes.  Nor do they indicate a range of standardized effect sizes that are not supported.

\subsection{Biases associated with confirmatory information processing}
To illustrate BFFs based on F statistics in replicated studies, we turn to a study reported in \cite{Fischer2008} that was replicated in 2015  \cite{OSF_article}.  Both studies sought to determine whether states of self-regulation depletion or ego threat caused participants to exhibit more bias in confirmatory information processing. The studies compared preferences for decision-consistent and decision-inconsistent information processing between three groups: high depletion of self-regulation, low depletion of self-regulation, and ego-threatened subjects.  The dependent variable consisted of a normalized score for participants' selection of decision-consistent and decision-inconsistent reports regarding a hiring decision upon which they had made a preliminary decision.  The original study's authors recruited 85 undergraduate students as subjects, while 140 subjects participated in the replicated study.  Differences between the outcomes in the three groups were assessed using one-way ANOVA.  The $F$ statistics reported in the original and replicated studies were $F_{2,82} = 4.05$ and $F_{2,137} =1.99$, respectively.

To construct a Bayes factor from these studies, it is first necessary to define the hypotheses being tested.  To make the discussion more general, we assume a total of $S$ studies; $S=2$ in this example.

Let $x_1,\dots,x_S$ denote $S$ independent $F$ statistics with numerator degrees-of-freedom $k_1,\dots,k_S$ and denominator degrees-of-freedom $m_1,\dots,m_S$, respectively.  In the present case, $k_1=k_2=2$, and $m_1 = n_1 -3$ and $m_2=n_2-3$, where $n_1=85$ and $n_2=140$ are the sample sizes in the two studies. 

Under the null hypothesis, we assume 
\begin{equation}
H_0: x_s \sim F(k_s,m_s,0), \qquad \mbox { for } s=1,\dots,S.
\end{equation}
Given the independence of the $\{ x_s \}$, the marginal density of the data under the null hypothesis is
\begin{equation}\label{mf0}
m_0(x_1,\dots,x_S) = \prod_{s=1}^S m_0(x_s) = \prod_{s=1}^S f(x_s \con k_s, m_s, 0).
\end{equation}

Under the alternative hypothesis, we assume
\begin{equation}
H_1: x_s \con \lambda \sim F(k_s,m_s,\lambda_s), \qquad \lambda_s \con \tau_{s}^2 \sim G\left(\frac{k_s}{2}+1,\frac{1}{2\tau_{s}^2} \right), \qquad \tau_s^2  = \frac{n_s \tilde{\omega}^2}{2} .
\end{equation}
Different prior distributions are specified for the non-centrality parameters $\{\lambda_s\}$ to account for the dependence of these parameters on each study's sample size.  However, the rate parameters $\tau^2_s$ that define these prior densities were determined from a common RMSES, $\tilde{\omega}$.  This stipulation models the belief that the  interventions have similar effects across studies.

Assuming that the non-centrality parameters are conditionally independent across studies, it follows that the marginal density of the data under the alternative hypothesis is
 \begin{equation}\label{mf1}
m_1(x_1,\dots,x_S \con \tilde{\omega} ) = \prod_{s=1}^S m_1(x_s \con \tau^2_{\tilde{\omega}}).
\end{equation}
Here, the dependence of the marginal densities on the assumed value of $\tilde{\omega}$ and $\tau^2_{\tilde{\omega}}$ has been indicated.   
Dividing equation (\ref{mf1}) by (\ref{mf0}) leads to
\begin{equation}\label{meta}
BF_{10}(x_1,\dots,x_S \con \tilde{\omega}) = \frac{\prod_{s=1}^S m_1(x_s \con \tau^2_{\tilde{\omega}})}{ \prod_{s=1}^S m_0(x_s)} 
= \prod_{s=1}^S \frac{m_1(x_s \con \tau^2_{\tilde{\omega}}  )}{m_0(x_s)} 
=  \prod_{s=1}^S BF_{10}(x_s \con \tau^2_{\tilde{\omega}}).
\end{equation} 
Thus, the Bayes factor for the combined study can be obtained by multiplying the Bayes factors from the individual studies.  The BFF is constructed by calculating the Bayes factors over a range of $\tilde{\omega}$.

Equation~(\ref{meta}) can be applied generally to obtain Bayes factors based on independent $z$, $t$, $\chi^2$, and $F$ statistics using Theorems 1--4 and Table 1, under the assumption that non-centrality parameters are drawn independently from their prior distributions.

\begin{figure}[htb] 
   \centering
\includegraphics[width=4.5in]{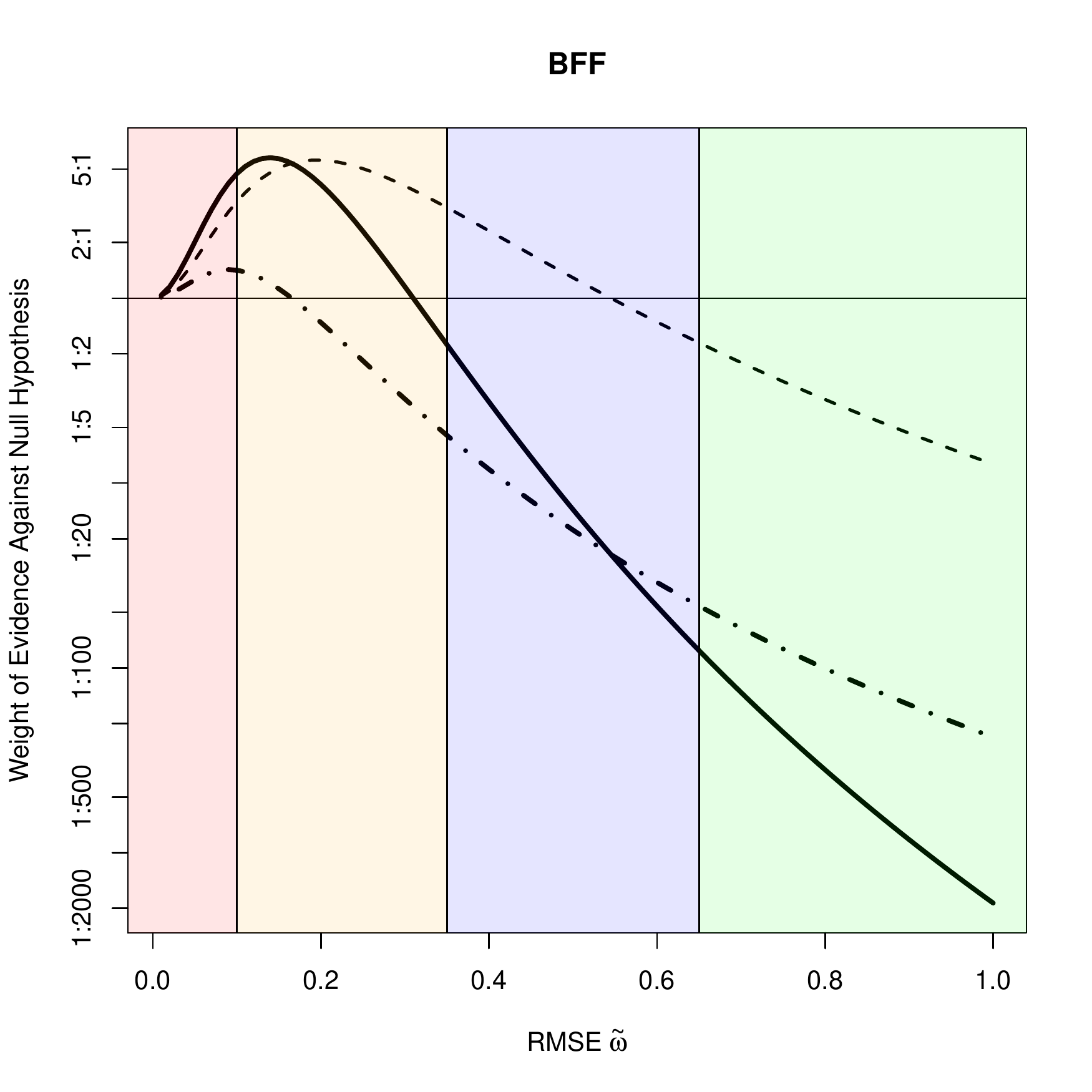}
   \caption{BFFs for confirmatory information processing studies.  The upper dotted line depicts the BFF from the original study; the lower dotted line for the replicated study.  The solid line represents the BFF obtained by multiplying the Bayes factors from the two studies.}
   \label{fig:F}
\end{figure}

Returning to our example, the Bayes factors for the two studies can be combined according to equation (24) and Theorem~4 as follows:
\begin{equation}
\begin{split}
 BF_{10}[f_{2,82}&=4.05,f_{2,137}=1.99\con \tilde{\omega}] = \\
 & BF_{10}\left[f_{2,82}=4.05\con \tau^2_{\tilde{\omega}} = \frac{85 \tilde{\omega}^2}{2}\right] \times  BF_{10}\left[f_{2,137}=1.99 \con \tau^2_{\tilde{\omega}} = \frac{140 \tilde{\omega}^2}{2}\right] .
 \end{split}
\end{equation}
Fig.~\ref{fig:F} depicts the BFF versus RMSES $\tilde{\omega}$ based on both experiments alongside the original and replication studies.  By combining information from the two studies, we see that there is support for very small or small differences in standardized effect sizes, with greater than 2:1 support for RMSES greater than 0.05 and less than 0.28, but no support for standardized effect sizes greater than 0.28.  The maximum Bayes factor against the null hypothesis of no effect was obtained for an RMSES of 0.14, where it was 5:1 against the null hypothesis of no group effect.

\section{Discussion}
The Bayes factors and BFFs described above are based on the specification of normal-moment and gamma prior densities imposed on scalar non-centrality parameters. Other prior specifications on non-centrality parameters are, of course, possible.  However, the proposed prior densities possess several attractive features.  Among these, they represent NAP densities, making it possible to accumulate evidence more rapidly in favor of true null hypotheses. They also yield closed-form expressions for Bayes factors, which facilitates BFF calculation. The coefficients-of-variation of the gamma priors in Theorems 3 and 4 are equal to $\sqrt{2/(k+2)}$ and so depend only on the (numerator) degrees of freedom of the test statistics.  
Under the proposed framework, the standard deviations of prior distributions on non-centrality parameters is thus scaled according to sample size.  
When used in conjunction with the normal moment priors specified in Theorems 1 and 2, these choices also yield Bayes factors that are invariant to the choice of test statistic in the sense that $z$ and $z^2=\chi_1^2$ tests, and $t_\nu$ and $t^2_\nu = F_{1,\nu}$ tests, produce the same Bayes factors when a common value of $\tau^2$ is selected.

For the test statistics considered above, it is possible to compute a ``maximum BFF" as the ratio of that test statistic's non-central alternative density to its central density under the null hypothesis (without averaging over a prior density).  This procedure essentially produces a plot of the likelihood ratio for each test statistic.  Because the probability that the test statistic matches this maximum value is either zero (continuous data) or small (discrete data), the maximum Bayes factor reported from such a procedure overstates evidence in favor of the alternative hypothesis.  This procedure also precludes the collection of evidence in favor of true null hypotheses and fails to model variability of standardized effect sizes across studies.

Many scientists now acknowledge the critical role that replication studies play in improving the reproducibility of scientific studies \cite{Munafo2017,Ioannidis2018}.  The final example demonstrates that BFFs provide a formal mechanism to combine information collected across replicated experiments using only reported test statistics. Their use thus provides a potential tool for enhancing the reproducibility of scientific research.  

An R package to calculate default BFFs for tests described in this article, ``BFF,'' is available for download at {\em cran.r-project.org}.

%\nocite{Abramowitz1970,Bishop1975,Mitra1958,Kendall1991}

\bibliography{article_ref}
\bibliographystyle{Science}


\section*{Acknowledgments} We thank M. Pourahmadi and A. Bhattarcharya for careful reading of the manuscript and helpful comments.  This research was supported by NIH CA R01 158113.  %{\bf Author contributions:} V.E.J. conceptualized the research and drafted the manuscript.  V.E.J. and S.P. proved the theorems.  R.S. created implementation software and identified applications.  All authors participated in editing and revising the manuscript.}

%\newpage

\section{Supplemental Material}

\subsection{Proofs of Theorems}
Throughout, $m_0(\bfx)$ and $m_1(\bfx)$ denote marginal densities of data vectors $\bfx$ under the null and alternative hypotheses, respectively, $BF_{10}(\bfx) = m_1(\bfx)/m_0(\bfx)$, and $\tau^2$ is treated as a constant.  The proofs of the first two theorems are similar to the proofs provided in the supplemental material of (24). 

\bigskip

\noindent{\em Proof of Theorem 1.}
By assumption, 
\begin{equation}
m_0(z) = c\exp\left(-z^2/2\right),  \qquad \mbox{where} \qquad c=\frac{1}{\sqrt{2\pi}}
\end{equation}
and 
\begin{eqnarray}
m_1(z\con\tau^2) &=& \int_{-\infty}^{\infty} \frac{c}{\sqrt{2\pi} \tau^3} \lambda^2 \exp\left( -\frac{\lambda^2}{2\tau^2} \right) \, \exp\left(-\frac{(z-\lambda)^2}{2}\right) d\lambda \\
&=& \int_{-\infty}^{\infty} \frac{c}{\sqrt{2\pi} \tau^3 } \lambda^2 
\exp\left[ -\frac{1}{2}\left( \frac{\lambda^2}{\tau^2} + (z-\lambda)^2 \right) \right] d\lambda .
\end{eqnarray}
Letting $a = \tau^2/(1+\tau^2)$ and noting 
\begin{equation}\label{onez_identity}
\frac{\lambda^2}{\tau^2} + (z-\lambda)^2 = \frac{1}{a} \left(
\lambda-az \right)^2 + \frac{az^2}{\tau^2} ,
\end{equation}
it follows that
\begin{eqnarray}
m_1(z\con \tau^2) &=&   \int_{-\infty}^{\infty} \frac{c}{\sqrt{2\pi} \tau^3 } \lambda^2 \exp\left\{ -\frac{1}{2} \left[ \frac{1}{a} \left(
\lambda-az \right)^2 +z^2 - a^2z^2 \right]\right\}d\lambda \\
&=& \frac{\sqrt{a} c}{\tau^3} \exp\left( -\frac{a z^2}{2\tau^2} \right) \int_{-\infty}^{\infty} \frac{1}{\sqrt{2\pi a}} \lambda^2 \exp\left[ -\frac{(\lambda-az )^2}{2a} \right] d\lambda .
\end{eqnarray}
The integral represents the second moment of a normal distribution with mean $az$ and variance $a$. Thus
\begin{equation}
m_1(z\con \tau^2) = \frac{\sqrt{a} c}{\tau^3}  \left[ a + (az)^2 \right]  \exp\left( \frac{a z^2}{2\tau^2}  \right).
\end{equation}
Substituting for $a$ and dividing $m_1(z \con \tau^2)$ by $m_0(z)$ produces the Bayes factor appearing in Theorem 1 of the article.
\begin{flushright}
$\blacksquare$
\end{flushright}

\bigskip
\noindent{\em Proof of Theorem 2.}
Under $H_0$, the marginal density of the $t$ statistic can be expressed
\begin{equation}\label{tdens}
m_0(t) = \frac{ \Gamma[(\nu+1)/2]}{\sqrt{\nu\pi} \Gamma(\nu/2)} \left( 1 + \frac{t^2}{\nu} \right)^{-(\nu+1)/2}.
\end{equation}
Defining
\begin{equation}
c = \frac{\nu^{\nu/2}}{\sqrt{\pi} \Gamma(\nu/2) d^{\nu+1} 2^{(\nu-1)/2}} \qquad \mbox{with} \qquad d = \sqrt{t^2+\nu}, 
\end{equation}
$m_1(t)$ can be expressed by integrating the integral form of the non-central $t$ density from (35) with respect to a $J(0,\tau^2)$ prior density to obtain
\begin{equation}\label{t1marg}
m_1(t) = \int_{-\infty}^{\infty} \int_0^{\infty} c \exp\left( -\frac{\nu\lambda^2}{2d^2}\right) y^{\nu} \exp\left[ \frac{(y-\lambda t/d)^2}{2} \right] \frac{\lambda^2}{\sqrt{2\pi} \tau^3} \exp\left(-\frac{\lambda^2}{2\tau^2}\right) dy \, d\lambda.
\end{equation}
Again letting $a=\tau^2/(\tau^2+1)$ and noting that
\begin{equation}
\frac{\nu\lambda^2}{d^2} + {(y-\lambda t/d)^2}+ \frac{\lambda^2}{\tau^2} = \frac{1}{a}\left( \lambda - \frac{ayt}{d}\right)^2 + y^2\left(1-\frac{at^2}{d^2}\right),
\end{equation}
application of Fubini's theorem implies that the integral with respect to $\lambda$ is proportional to the second moment of a normal density with variance $a$ and mean $(ayt/d)$.  Thus,
\begin{equation}\label{into}
m_1(t) = \int_0^{\infty} c y^{\nu} \frac{1}{\tau^2\sqrt{\tau^2+1}}\left[ a + \frac{(ayt)^2}{d^2}\right] \exp\left[-\frac{y^2}{2}\left(1-\frac{at^2}{d^2}\right) \right] dy.
\end{equation}
Changing variables to $z = y^2$, it follows that equation (\ref{into}) represents the weighted sum of two gamma functions.  Define
\begin{equation}
b = 1- \frac{at^2}{d^2}, \qquad \mbox{which implies} \qquad \frac{1}{b} = \frac{1 + \frac{t^2}{\nu}}{1+\frac{t^2}{\nu(\tau^2+1)}}.
\end{equation} 
Integrating equation (\ref{into}) leads to
\begin{eqnarray}
m_1(t) &=& c\, \frac{1}{2(\tau^2+1)^{3/2}} \left[ \frac{2^{(\nu+1)/2} \Gamma\left(\frac{\nu+1}{2}\right)}{b^{(\nu+1)/2}} +  
                 \frac{a t^2 2^{(\nu+3)/2}  \Gamma\left(\frac{\nu+3}{2}\right)}{ d^2 b^{(\nu+3)/2} }\right] \\
            &=& \frac{\nu^{\nu/2}}{\sqrt{\pi} \Gamma(\nu/2) d^{\nu+1} (\tau^2+1)^{3/2}} \left[  \frac{\Gamma\left(\frac{\nu+1}{2}\right)}{b^{(\nu+1)/2}} + 
                   \frac{at^2(\nu+1)\Gamma\left(\frac{\nu+1}{2}\right)}{d^2 b^{(\nu+3)/2}} \right] .\label{intt}
\end{eqnarray}
Also,
\begin{equation}
\left(d \nu^{-1/2}\right)^{-\nu-1} = \left( 1+\frac{t^2}{\nu}\right)^{-(\nu+1)/2} .
\end{equation}
Thus, dividing $m_1(t)$ in equation (\ref{intt}) by $m_0(t)$ in equation (\ref{tdens}) leads to
\begin{equation}
BF_{10}(t) = (\tau^2+1)^{-3/2} b^{-(\nu+1)/2} \left[ 1 + \frac{a(\nu+1)t^2}{bd^2} \right].
\end{equation} 
But $1/b = r/s$ and 
\begin{equation}
\frac{a(\nu+1)t^2}{bd^2} = \frac{qt^2}{s},
\end{equation}
yielding the Bayes factor specified in the theorem.

\begin{flushright}
$\blacksquare$
\end{flushright}

\bigskip
\noindent{\em Proof of Theorem 3.}
Under $H_0$, the density of $h$ is
\begin{equation}
m_0(h) = \frac{1}{2^{k/2}\, \Gamma\left(\frac{k}{2}\right)} h^{\frac{k}{2}-1} \exp\left(-\frac{h}{2}\right).
\end{equation}

Under $H_1$, the marginal density of $h$ is obtained by integrating the non-central chi-squared density (35) with respect to the prior on the non-centrality parameter $\lambda$:
\begin{eqnarray}
m_1(h) &=& \int_0^\infty \sum_{i=0}^\infty \frac{e^{-\lambda/2} \left(\frac{\lambda}{2}\right)^i}{i! \, 2^{\frac{k+2i}{2}} \Gamma\left(\kt+i\right)} h^{\kt+i-1} \exp\left(-\frac{h}{2}\right) \\
& & \ \ \ \ \times \frac{\lambda^{\frac{k}{2}} \exp\left(-\frac{\lambda}{2\tau^2}\right)}{(2\tau^2)^{\kt+1}\Gamma\left(\frac{k}{2}+1\right)} \ d\lambda.
\end{eqnarray}
Noting that
\begin{equation}
\int_0^\infty \lambda^{\kt+i} \exp\left[ -\lambda\left(\frac{1}{2}+\frac{1}{2\tau^2}\right) \right] d\lambda = 
 \frac{ \left(\kt+i \right)\Gamma\left(\kt+i\right)}{\left( \frac{{\tau^2+1}}{2\tau^2}\right)^{\kt+i+1}} ,
\end{equation}
application of Fubini's theorem leads to 
\begin{equation}\label{sery}
m_1(h) = b(h) \sum_{i=0}^\infty \left( \frac{2\tau^2}{{\tau^2+1}}\right)^i  \frac{(\kt+i)h^i}{2^{2i}\, i!},
\end{equation}
where
\begin{equation}
b(h) = \frac{1}{({\tau^2+1})^{\kt+1} \Gamma\left(\kt+1\right) 2^\kt}  \ h^{\kt-1} e^{-h/2}. 
\end{equation}
Because
\begin{equation} 
\sum_{j=0}^\infty \frac{(c+j) a^j}{j!} = (c+a)e^a,
\end{equation}
equation (\ref{sery}) can be rewritten as
\begin{equation}\label{three}
m_1(h) = b(h) \left[ \kt  + \frac{\tau^2 h}{2({\tau^2+1})} \right] \exp\left(\frac{\tau^2 h}{2({\tau^2+1})}\right).
\end{equation}
Noting that 
\begin{equation}
\frac{b(h)}{m_0(h)} = ({\tau^2+1})^{-(\kt+1)} \frac{2}{k}
\end{equation}
and dividing in equation (\ref{three}) yields the Bayes factor stated in the theorem.
\begin{flushright}
$\blacksquare$
\end{flushright}
As an aside, the $m_1(h)$ can be expressed as
\begin{equation}
m_1(h) = \frac{1}{{\tau^2+1}} g\left(h\,; \kt, \frac{1}{2({\tau^2+1})}\right) + \frac{\tau^2}{{\tau^2+1}} g\left(h\,; \kt+1, \frac{1}{2({\tau^2+1})}\right).
\end{equation}

\bigskip

\noindent{\em Proof of Theorem 4.}
Under $H_0$, the marginal density of the $f$ statistic is 
\begin{equation}
m_0(f) = B\left(\frac{k}{2},\frac{m}{2}\right)^{-1} \left(\frac{k}{m}\right)^{k/2} f^{k/2 -1 }\left( 1 + \frac{k}{m}f\right)^{-(k+m)/2},
\end{equation}
where $B(k/2,m/2)  = \Gamma(k/2)\Gamma(m/2)/\Gamma[(k+m)/2]$ is the beta function.
The marginal density of $f$ under $H_1$ is obtained by integrating the non-central $f$ density (35) with respect to the gamma prior on its non-centrality parameter $\lambda$:
\begin{eqnarray}
m_1(f) &=& \int_0^\infty  \left(\frac{k}{m}\right)^{k/2} e^{-\lambda/2} 
\sum_{r=0}^\infty \left(\frac{k\lambda}{2m}\right)^r \frac{1}{r!}
 B\left(\frac{k}{2}+r,\frac{m}{2}\right)^{-1} \\ 
 & & \times \ \ \ \ \frac{f^{r+k/2-1}}{\left(1+kf/m\right)^{r+(k+m)/2}} \ 
\frac{\lambda^{k/2}\exp(-\lambda/(2\tau^2))}{(2\tau^2)^{k/2+1}\Gamma(k/2)} \  d\lambda.
 \end{eqnarray}
 Defining 
 \begin{equation}
 a(f) = \frac{\left(\frac{k}{m}\right)^{k/2} f^{k/2 -1 } }{\Gamma(\frac{k}{2}+1)\Gamma(\frac{m}{2})(2\tau^2)^{1+k/2}\left( 1 + \frac{k}{m}f\right)^{-(k+m)/2}},
 \end{equation}
application of Fubini's theorem, recognizing
\begin{equation}
\Gamma\left(\frac{k}{2}+r\right) \left( \frac{1}{2} + \frac{1}{2\tau^2}\right)^{-r-k/2-1} = \int_0^{\infty} \lambda^{r+k/2} \exp\left[-\frac{\lambda}{2}\left(1+\frac{1}{\tau^2}\right) \right] d\lambda ,
\end{equation}
and expanding $\Gamma(r+\frac{k}{2}+\frac{m}{2})$ as an integral allows the marginal density to be expressed as
\begin{eqnarray}
m_1(f) &=&  \int_0^\infty a(f) \ e^{-x} \ \frac{x^{k/2+m/2-1}}{\left( \frac{1}{2}+\frac{1}{2\tau^2}\right)}  \\
& & \times \ \ \ \ \sum_{r=0}^\infty \left(r + \frac{k}{2}\right) 
\left[ \frac{kxf}{m\left(1+\frac{k}{m}f \right) \left( 1+\frac{1}{\tau^2} \right)} \right]^r \ \frac{1}{r!} \, dx \label{exp}.
\end{eqnarray}
The series in equation (\ref{exp}) can be summed and is equal to 
\begin{equation}
\frac{k}{2} e^{b(x)} + b(x)e^{b(x)},
\end{equation}
where $b(x)$ is the bracketed term raised to power $r$.  Letting $c=b(x)/x$, 
%\begin{equation}
%cx = B(x)\qquad \mbox{and} \qquad D=\left[\frac{2\tau^2}{\tau^2+1}\right]^{\frac{k}{2}},
%\end{equation}
the integration with respect to $x$ represents the sum of two scaled gamma functions, leading to
\begin{equation}
m_1(f) = \frac{a(f)\,\Gamma\left(\frac{k+m}{2}\right)}{(1-c)^{\frac{k+m}{2}}} \left[\frac{2\tau^2}{\tau^2+1}\right]^{\frac{k}{2}}\left[ \frac{k}{2} + \left(\frac{c}{1-c}\right)\frac{k+m}{2} \right].
\end{equation}
Noting that 
\begin{equation}
1-c = \left( \frac{1}{1+\frac{k}{m}f}\right) \left( \frac{1}{1+\frac{k}{m(1+\tau)}f}\right),
\end{equation}
further simplification leads to 
\begin{eqnarray}
m_1(f) &=& \left(\frac{k}{m}\right)^{k/2}  \frac{1}{(\tau^2+1)^{(\frac{k}{2}+1)}} \left[  \frac{1}{B\left( \frac{k}{2}, \frac{m}{2} \right)} 
\frac{f^{k/2-1}}{ \left[ 1+ \frac{kf}{m(1+\tau^2)} \right]^{\frac{k+m}{2}}} \right. \\
&+& % \left(\frac{k}{m}\right)^{k/2}   \frac{1}{(\tau^2+1)^{(\frac{k}{2}+1)}} 
\left.   \frac{k\tau^2}{m(\tau^2+1)} \frac{1}{B\left( \frac{k}{2}+1, \frac{m}{2} \right)} 
\frac{f^{k/2}}{ \left[ 1+ \frac{kf}{m(1+\tau^2)} \right]^{\frac{k+m}{2}+1}} \right].
\end{eqnarray}
Substituting $\frac{k}{k+m} B(\frac{k}{2},\frac{m}{2})$ for $B(\frac{k}{2}+1,\frac{m}{2})$ and dividing $m_1(f)$ by $m_0(f)$ yields the Bayes factor specified in theorem.

\begin{flushright}
$\blacksquare$
\end{flushright}

As an aside, if 
\begin{equation} 
\frac{Y}{1+\tau^2} \sim F(k,m,0) \qquad \mbox{and} \qquad \frac{kZ}{(k+2)(1+\tau^2)} \sim F(k+2,m,0),
\end{equation}
then 
\begin{equation}
m_1(f) =  \frac{1}{1+\tau^2} \ p_Y(f) + \frac{\tau^2}{1+\tau^2} \ p_Z(f).
\end{equation}

 \subsection{Selection of $\tau^2$}
 As stated in the main article, our criteria for selecting $\tau^2$ in the specification of the prior distribution on the non-centrality parameters of the test statistics is to select $\tau^2$ so that the prior modes on the non-centrality parameters correspond to a specified standardized effect size $\omega$.  For scalar variables, standardized effect sizes refer to the number of standard deviations the mean of a single observation falls from a specified null hypothesis, while for vector-valued parameters we refer to standardized effect size as the effect size standardized by the Cholesky decomposition of the (asymptotic) covariance matrix.  Rationale for the choice of $\tau_\omega^2$ for each of the tests listed in Table~1 follow.
 
 \subsubsection{One-sample $z$ and $t$ tests}
Suppose $x_i$, $i=1,\dots,n$, represent independent and identically distributed (iid) $N(\mu,\sigma^2)$ random variables.  Without loss of generality, the point null hypothesis may be specified as $H_0: \mu=0$ and either Theorem 1 ($\sigma^2$ known) or Theorem 2 ($\sigma^2$ unknown) may be applied, with test statistics $\sqrt{n}\bar{x}/\sigma$ or $\sqrt{n}\bar{x}/s$, respectively. In both cases, the non-centrality parameter is $\lambda = \sqrt{n}\mu/\sigma$. The modes of a $J(0,\tau^2)$ density are $\pm \sqrt{2} \tau$. Defining the standardized effect  
$\omega = \mu/\sigma$ and equating $\lambda^2 = 2\tau^2$ leads to $\tau_\omega^2 = n\omega^2/2$.

\subsubsection{Two-sample $z$ and $t$ tests}
Suppose $x_{ij}$, $i=1,\dots,n_j$, $j=1,2$, represent iid $N(\mu_j,\sigma^2)$ random variables. The point null hypothesis is $H_0: \mu_1=\mu_2$ and either Theorem 1 ($\sigma^2$ known) or Theorem 2 ($\sigma^2$ unknown) may be applied, with test statistics 
\begin{equation} 
z = \frac{\bar{x_1}-\bar{x_2}}{\sigma \sqrt{1/n_1+1/n_2}} \qquad \mbox{or} \qquad  t_{\nu} = \frac{\bar{x_1}-\bar{x_2}}{s \sqrt{1/n_1+1/n_2}},
\end{equation}
where $\nu=n_1+n_2-2$.

The standardized effect size for both tests is $\omega = (\mu_1-\mu_2)/\sigma$, and the non-centrality parameter is
\begin{equation}
\lambda = \frac{\mu_1-\mu_2}{\sigma\sqrt{1/n_1+1/n_2}} = \frac{\sqrt{n_1n_2}\omega}{\sqrt{n_1+n_2}},
\end{equation}
Equating $\lambda = \sqrt{2}\tau $ leads to 
\begin{equation}
\tau^2 = \frac{ n_1n_2\, \omega^2}{2(n_1+n_2)}.
\end{equation}

\subsubsection{Multinomial and Poisson Tests}
We assume that the asymptotic conditions described in, for example, (36) are satisfied for chi-squared test statistics of the form
\begin{equation}
h(\hat{\bt}) = \sum_{k=1}^{K} \frac{[n_k - n f_k(\hat{\bt})]^2}{n f_k(\hat{\bt})},
\end{equation} 
defined for counts $n_k$ in $K$ cells, an $s$ dimensional parameter $\theta \in \Theta$ with $1\leq s < K$, and a $k\times 1$ vector-valued function ${\bf f}(\bt) = (f_1(\bt),\dots,f_K(\bt))'$ that satisfies $\sum_k f_k(\bt) = 1$.  Given $n = \sum_{k=1}^K  n_k$, the counts $\{ n_k\}$ are assumed to follow a multinomial distribution.  In particular, this assumption is satisfied if the counts $n_k$ arise as independent Poisson random variables.   Under the null hypothesis, the cell probabilities satisfy $H_0: \bpi = f_k(\bt)$ for some $\bt \in \Theta$.  Under the alternative hypothesis, the true cell probabilities ${\bf p}$ are assumed to satisfy ${\bf p}-\bpi = O(n^{-1/2})$.  If $\hat{\bt}$ is the maximum likelihood or minimum chi-squared estimate of $\bt$, then under certain regularity conditions, (37) showed that $h(\hat{\bt})$ converges to a non-central chi-squared distribution on $K-s-1$ degrees of freedom and non-centrality parameter 
\begin{equation}
n (\bp - \bpi)' {\bf D}_\bpi^{-1} (\bp-\bpi) , \quad \mbox{where} \quad {\bf D}_\bpi = \{ \mbox{diag}(\bpi) \}. 
\end{equation}
We define the standardized effect size vector as 
\begin{equation}
\bfo = \left\{ \frac{p_k - \pi_k}{\sqrt{\pi_k}} \right\}.
\end{equation}
The mode of the prior density on the non-centrality parameter in Theorem~3 occurs at $k\tau^2$, and the non-centrality parameter can be written as $n\bfo'\bfo$.  Matching the mode of the prior density to the non-centrality parameter leads to 
\begin{equation}
\tau^2 = \frac{n\bfo'\bfo}{k}.
\end{equation}
In some applications, it is convenient to replace $\bfo'\bfo$ by $k\tilde{\omega}^2$, where $\tilde{\omega}^2 = \frac{1}{k} \sum_{i=1}^k \omega_i^2$, the average squared standardized effect.     

\subsubsection{Linear Models}
We assume that the $n\times 1$ data vector ${\bf y}$ satisfies
\begin{equation}
{\bf y} \sim N({\bf X} \bfbeta, \sigma^2 {\bf I}_n)
\end{equation}
for a $p\times 1$ vector $\bfbeta$ and $n\times p$ matrix ${\bf X}$ of rank $p<n$.
The null hypothesis can be expressed as $H_0: {\bf A}\bfbeta  = {\bf a}$ where the rank of ${\bf A}$ is $k\leq p$.  The $F$ statistic against the alternative hypothesis $H_1: {\bf A}\bfbeta  \neq {\bf a}$ can be expressed as
\begin{equation} 
F_{k,n-p} = \frac{(RSS_0-RSS_1)/k}{RSS_1/(n-p)},
\end{equation}
where $RSS_0$ is the constrained residual sum-of-squares under the null model and $RSS_1$ is the residual sum-of-squares under the unconstrained model.

As demonstrated in, for example (25), the non-centrality parameter for the $F$ statistic can be expressed as
\begin{equation} 
\lambda = \frac{n({\bf A}\bfbeta - {\bf a})'{\bf V}^{-1} ({\bf A}\bfbeta - {\bf a})}{2\sigma^2}, \qquad \mbox{where} \quad {\bf V} = \left[{\bf A}\left(\frac{{\bf X}'{\bf X}}{n}\right)^{-1}{\bf A}'\right].
\end{equation}
Letting ${\bf L}$ denote the Cholesky decomposition of ${\bf V}$, i.e., ${\bf L}{\bf L}' = {\bf V}$, we define the standardized effect size vector $\bfo$ as
\begin{equation} 
\bfo = \frac{{\bf L}^{-1}({\bf A}\bfbeta - {\bf a})}{\sigma}.
\end{equation}
It follows that the non-centrality parameter can be written
\begin{equation}
\lambda = \frac{n\bfo' \bfo}{2}.
\end{equation}
Equating $\lambda= k\tau^2$ leads to  
\begin{equation}
\tau^2= \frac{n\bfo' \bfo}{2k}.
\end{equation}

As for count models, it is sometimes convenient to replace $\bfo'\bfo$ by $k\tilde{\omega}^2$.

\subsubsection{Likelihood Ratio Statistic}
We assume that the $n\times 1$ data vector ${\bf y}$ is generated from a parametric family of densities indexed by $\bt = (\bt_r,\bt_s)$, say $f({\bf y},\bt)$.  The likelihood function is denoted by $l(\bt \con {\bf y})$.  The null hypothesis is $H_0: \bt_r = \bt_{r0}$, where $\dim(\bt_r) = k$. The likelihood ratio statistic is 
\begin{equation}
h = -2 \log\left[ \frac{l(\bt_{r0},\hat{\bt_s})}{l(\hat{\bt})} \right], 
\end{equation}
where $\hat{\bt}$ is the unconstrained maximum likelihood estimate (MLE) and $\hat{\bt_s}$ is the constrained MLE under $H_0$.

Assuming regularity conditions in (38), $h$ converges to a $\chi^2$ distribution on $k$ degrees of freedom and non-centrality parameter
\begin{equation}
\lambda =   n(\bt_r - \bt_{r0})' {\bf V}_r^{-1} (\bt_r - \bt_{r0}) \qquad \mbox{where} \qquad 
{V}_{r_i,r_j}^{-1} = - \frac{1}{n}E\left[ \frac{ \partial^2 \log l}{\partial \theta_{r_i} \partial \theta_{r_j} } \right].
\end{equation}
Let ${\bf L}$ denote the Cholesky decomposition of ${\bf V}_r$  and define
\begin{equation}
\bfo = {\bf L}^{-1}(\bt_r-\bt_{r0}).
\end{equation}
Then $\lambda = n\bfo'\bfo$, and setting $\lambda = k\tau^2$ implies
\begin{equation}
\tau^2 = \frac{n\bfo'\bfo}{k}.
\end{equation}
As for the linear and discrete models, $\bfo'\bfo$ may be replaced by $k \tilde{\omega}^2$.

\end{document}